\documentclass[a4paper]{article}

\usepackage[english]{babel}
\usepackage[utf8x]{inputenc}
\usepackage{enumerate,amsmath,amssymb,amstext,amsthm}
\usepackage{graphicx}
\usepackage[colorinlistoftodos]{todonotes}
\usepackage{authblk}

\usepackage{hyperref}

\newtheorem{theorem}{Theorem}[section]
\newtheorem{lemma}[theorem]{Lemma}
\newtheorem{corollary}[theorem]{Corollary}

\theoremstyle{definition}
\newtheorem*{definition}{Definition}

\newcommand{\bartheta}{\bar{\vartheta}}

\renewcommand{\l}{\lambda}
\newcommand{\m}{\mu}
\renewcommand{\th}{\theta}
\newcommand{\btau}{\bar{\tau}}
\newcommand{\bth}{\bar{\th}}
\newcommand{\bm}{\bar{\m}}
\newcommand{\bl}{\bar{\l}}
\newcommand{\bk}{\bar{k}}

\newcommand{\qeds}{\qed\vspace{.2cm}}
\DeclareMathOperator{\tr}{Tr}

\DeclareMathOperator{\im}{Im}

\title{Homomorphisms of Strongly Regular Graphs}
\author{David E.~Roberson}
\affil{Department of Computer Science \\ University College London}
\date{}

\begin{document}
\maketitle
\begin{abstract}
We prove that if $G$ and $H$ are primitive strongly regular graphs with the same parameters and $\varphi$ is a homomorphism from $G$ to $H$, then $\varphi$ is either an isomorphism or a coloring (homomorphism to a complete subgraph). Moreover, any such coloring is optimal for $G$ and its image is a maximum clique of $H$. Therefore, the only endomorphisms of a primitive strongly regular graph are automorphisms or colorings. This confirms and strengthens a conjecture of Peter Cameron and Priscila Kazanidis that all strongly regular graphs are cores or have complete cores. The proof of the result is elementary, mainly relying on linear algebraic techniques. In the second half of the paper we discuss implications of the result and the idea underlying the proof. We also show that essentially the same proof can be used to obtain a more general statement.

\end{abstract}

\section{Introduction}\label{sec:intro}

A \emph{homomorphism} between two graphs $G$ and $H$ is a function $\varphi: V(G) \to V(H)$ such that $\varphi(u) \sim \varphi(v)$ whenever $u \sim v$, where `$\sim$' denotes adjacency. Whenever a homomorphism exists from $G$ to $H$, we write $G \to H$, and if both $G \to H$ and $H \to G$ then we say that $G$ and $H$ are \emph{homomorphically equivalent}. Given a homomorphism $\varphi$ from $G$ to $H$, we will abuse terminology somewhat and refer to the subgraph of $H$ induced by $\{\varphi(u) : u \in V(G)\}$ as the \emph{image} of $\varphi$, and denote this by $\im \varphi$. It is easy to see that a $c$-coloring of a graph $G$ is equivalent to a homomorphism from $G$ to the complete graph on $c$ vertices, $K_c$. More generally, we will refer to any homomorphism whose image is a clique (complete subgraph) as a \emph{coloring}.

A homomorphism from a graph $G$ to itself is called an endomorphism, and it is said to be \emph{proper} if it is not an automorphism of $G$, or equivalently, its image is a proper subgraph of $G$. A graph with no proper endomorphisms is said to be a \emph{core}, and these play a fundamental role in the theory of homomorphisms since every graph is homomorphically equivalent to a unique core. We refer to the unique core homomorphically equivalent to $G$ as \emph{the core of $G$}. It is known~\cite{AGT}, and not difficult to show, that the core of $G$ is isomorphic to any vertex minimal induced subgraph of $G$ to which $G$ admits an endomorphism.

If the core of a graph $G$ is a complete graph $K_c$, then $G$ must contain a clique of size $c$ and must also be $c$-colorable. Therefore, $\omega(G) = \chi(G) = c$. Conversely, if $\omega(G) = \chi(G) = c$, then the core of $G$ is $K_c$. If a graph is either a core or has a complete graph as a core, then it is said to be \emph{core-complete}. Many known results on cores are statements saying that all graphs in a certain class are core-complete~\cite{symcores, geocores, altcores}, and often it remains difficult to determine whether a given graph in the class is a core or has a complete core.

For some classes of graphs, something stronger than core-completeness can be shown. A graph $G$ is a \emph{pseudocore} if every proper endomorphism of $G$ is a coloring. It follows that such a graph either has no proper endomorphisms and is thus a core, or has some proper endomorphism to a clique and thus has a complete core. In other words, any pseudocore is core-complete, although the converse does not hold (consider a complete multipartite graph). Similarly, it is easy to see that any core is a pseudocore, but the converse does not hold in this case either (for instance the Cartesian product of two complete graphs of equal size at least three).

In this paper, we will focus on homomorphisms and cores of strongly regular graphs. An $n$-vertex $k$-regular graph is said to be \emph{strongly regular} with parameters $(n,k,\l,\m)$ if every pair of adjacent vertices has $\l$ common neighbors, and every pair of distinct non-adjacent vertices has $\m$ common neighbors. For short, we will call such a graph an $SRG(n,k,\l,\m)$. A strongly regular graph is called \emph{imprimitive} if either it or its complement is disconnected. In such a case, the graph or its complement is a disjoint union of equal sized complete graphs. Homomorphisms of these graphs are straightforward, and so we will only consider \emph{primitive} strongly regular graphs here. Because of this, from now on when we consider a strongly regular graph, we will implicitly assume that it is primitive. In this case, we always have that $1 \le \m < k$, and that the diameter is two.

Cameron \& Kazanidis~\cite{symcores} showed that a special class of strongly regular graphs, known as rank 3 graphs, are all core-complete. A graph is \emph{rank 3} if its automorphism group acts transitively on vertices, ordered pairs of adjacent vertices, and ordered pairs of distinct non-adjacent vertices. The rank refers to the number of orbits on ordered pairs of vertices, and so after complete or empty graphs, rank 3 graphs are in a sense the graphs with the most symmetry. The proof of Cameron \& Kazanidis exploits this symmetry by noting that either no pair of non-adjacent vertices can be identified (mapped to the same vertex) by an endomorphism of a rank 3 graph, or every such pair can. In the former case, the graph must be a core. In the latter, any endomorphic image that contains non-adjacent vertices cannot be minimal, and therefore the core must be complete.

Strongly regular graphs can be viewed as combinatorial relaxations of rank 3 graphs and, following their result, Cameron \& Kazanidis (tentatively) conjectured that all strongly regular graphs are core-complete. Towards this, Godsil \& Royle~\cite{geocores} showed that many strongly regular graphs constructed from partial geometries are core-complete. A partial geometry is simply a point-line incidence structure obeying certain rules. The point graph of a partial geometry has the points as vertices, such that two are adjacent if they are incident to a common line. The properties of partial geometries guarantee that their point graphs are strongly regular, and they are typically referred to as \emph{geometric} graphs.

Godsil \& Royle showed that the point graphs of generalized quadrangles are pseudocores, as are the block graphs of $2$-$(v,k,1)$ designs and orthogonal arrays with sufficiently many points. As they note, a result of Neumaier~\cite{neumaier} is that for a fixed least eigenvalue, all but finitely many strongly regular graphs are the block graphs of $2$-$(v,k,1)$ designs or orthogonal arrays. Thus their result makes a significant step towards the conjecture of Cameron \& Kazanidis. The main idea used in the proof of the Godsil \& Royle result is that any endomorphism must map maximum cliques to maximum cliques. Starting with this simple observation, they show that if $G$ is geometric, and the maximum cliques of $G$ are exactly the lines of the underlying partial geometry, then $G$ is a pseudocore. It then remains to show when this assumption on the maximum cliques holds true.



The main result of this paper is that if $G$ and $H$ are both strongly regular graphs with parameters $(n,k,\l,\m)$, and $\varphi$ is a homomorphism from $G$ to $H$, then $\varphi$ is either an isomorphism or a coloring. Letting $G = H$, this statement implies that all strongly regular graphs are pseudocores, thus proving and strengthening the conjecture of Cameron \& Kazanidis. Using our main result and some previously known results, we also show that in the case where $\varphi$ is a coloring, we must have $\chi(G) = \omega(H)$ and this value is equal to the Hoffman bound on chromatic number which depends only on $(n,k,\l,\m)$. It follows from this that any strongly regular graph $G$ falls into one of four classes depending what subset of $\{\omega(G), \chi(G)\}$ meets the Hoffman bound. Using this we show that the homomorphism order of strongly regular graphs with a fixed parameter set has a simple description.

We also prove a generalization of our main result, where the strong regularity assumption on $H$ is replaced by a strictly weaker algebraic condition. In this more general case, we are only able to conclude that any homomorphism from $G$ to $H$ is either a coloring or an isomorphism to an induced subgraph of $H$.


The original idea and the inspiration for the proof of the main result comes from the theory of vector colorings, which are a homomorphism-based formulation of the famous Lov\'{a}sz theta function. The author was aided greatly by a collaboration with Chris Godsil, Brendan Rooney, Robert \v{S}amal, and Antonios Varvitsiotis which produced three papers~\cite{UVC1, UVC2, UVC3} on vector colorings. In particular, the second paper~\cite{UVC2} focused specifically on using vector colorings to restrict the possible homomorphisms between graphs. Note however that we will present an elementary proof of our main result which only requires basic knowledge of linear algebra and certain aspects of strongly regular graphs which we will review in Section~\ref{sec:SRGprops}. The connection between the proof techniques and vector colorings will not be discussed until Section~\ref{sec:vectcols}.

Although the main concrete contribution of this paper is the resolution and strengthening of the Cameron \& Kazanidis conjecture, we believe that the real significance of this work is the step it takes towards understanding how combinatorial regularity can impact the endomorphisms and core of a graph. Symmetry conditions, such as vertex- or distance-transitivity, often have easy-to-derive consequences for the endomorphisms and/or core of a graph. This is perhaps not surprising, since such symmetry conditions are assumptions about the automorphisms of a graph, which are just special cases of endomorphisms. However, it appears to be more challenging to make use of analogous regularity conditions, such as being strongly or distance regular. In fact, we believe that ours is the first example of such a result. Interestingly, by showing that strongly regular graphs are pseudocores, we establish a stronger result than was previously known even under the more stringent symmetry condition of being rank 3. Moreover, we know of no way to directly use the assumption of being rank 3 to show that a graph is a pseudocore.

%

\subsection{Notation}

Whenever we use $\th$ and $\tau$, we will be referring to the second largest and minimum eigenvalues of a strongly regular graph. This will sometimes be done without explicitly stating it. We will also use $m_\th$ and $m_\tau$ to denote the multiplicities of these eigenvalues, and $E_\th$ and $E_\tau$ will refer to the projections onto the corresponding eigenspaces.

The all ones matrix will be denoted by $J$. For a matrix $M$, we will use $\text{sum}(M)$ to refer to the sum of the entries of $M$. For two matrices $M$ and $N$ with the same dimensions, $M \circ N$ will denote their Schur, or entrywise, product.

The complement of a graph $G$ will be denoted by $\overline{G}$, and more generally we will add a bar over usual notation to refer to the analog in the complement. For instance, $\bth$ will refer to the second largest eigenvalue of the complement of a given strongly regular graph.

As already noted, we will use $u \sim v$ to mean that $u$ and $v$ are adjacent vertices. We will also use $u \not\sim v$ when $u$ and $v$ are not adjacent, which includes the case where $u = v$ since a vertex is not adjacent to itself. Sometimes we will need to exclude the $u = v$ case, and for this we will use $u \not\simeq v$. We will also refer to $u$ and $v$ as non-neighbors whenever $u \not\simeq v$. Lastly, note that $u \not\simeq v$ is equivalent to $u$ and $v$ being adjacent in the complement graph.

\section{Properties of Strongly Regular Graphs}\label{sec:SRGprops}

Here we will introduce some basic properties of strongly regular graphs that we will need later. We do not aim to give a full proof of every result, but rather enough explanation for the interested reader to work out the details. Most of these results are standard, and can be found in~\cite{AGT} or even on some widely used online sources that are not considered citable. Those familiar with strongly regular graphs can probably skip this section, with the possible exception of Lemma~\ref{lem:N2connected} and the definition of the cosines of a strongly regular graph at the end of Section~\ref{subsec:projections}.

\subsection{Algebraic properties}\label{subsec:algprops}

Let $G$ be an $SRG(n,k,\l,\m)$ with adjacency matrix $A$. Since $G$ is a connected $k$-regular graph, $k$ is a simple eigenvalue of $A$ with the all-ones vector as its unique (up to scalar) eigenvector. In particular, this implies that $AJ = JA = kJ$, 
and $Jz = 0$ for any $z$ that is an eigenvector for any eigenvalue of $A$ other than $k$.

As with any graph, the $uv$-entry of $A^m$ counts the number of walks of length $m$ between vertices $u$ and $v$. Using this and the definition of strongly regular graphs, it is not difficult to see that the matrix $A$ must satisfy the following:
\begin{equation}\label{eqn:evals}
A^2 + (\m - \l)A + (\m - k)I = \m J.
\end{equation}
Multiplying both sides of the above by an eigenvector of $A$ for an eigenvalue other than $k$, we see that all eigenvalues of $A$ other than $k$ must satisfy the equation $x^2 + (\m - \l)x + (\m - k) = 0$. Since $\tr(A) = 0$, the sums of the eigenvalues of $A$ must be zero, and this can be used to show that both roots of the above polynomial do occur as eigenvalues of $A$. Therefore, the other two distinct eigenvalues of an $SRG(n,k,\l,\m)$, denoted $\th$ and $\tau$, depend only on the parameters and are given as follows:
\begin{align*}
\th &= \frac{1}{2}\left[(\l - \m) + \sqrt{(\l - \m)^2 + 4(k - \m)}\right] \\
\tau &= \frac{1}{2}\left[(\l - \m) - \sqrt{(\l - \m)^2 + 4(k - \m)}\right].
\end{align*}
Note that these eigenvalues satisfy $k > \th > 0 > \tau$. The multiplicities, $m_\th$ and $m_\tau$, of $\th$ and $\tau$ can also be expressed in terms of the parameters $n,k,\l,\m$, but we will not need to make their values explicit. A key point to take away from this is that the eigenvalues, including their multiplicities, of a strongly regular graph depend only on the parameters, not on the specific graph.


\subsection{Projections onto eigenspaces}\label{subsec:projections}

For any real symmetric matrix $M$ with distinct eigenvalues $\zeta_1, \ldots, \zeta_m$, the projector onto the $\zeta_i$-eigenspace is a polynomial in $M$. This can be easily seen by considering how the matrix
\[\prod_{j \ne i} \frac{1}{\zeta_i - \zeta_j}(M - \zeta_j I)\]
acts on an orthogonal basis of eigenvectors of $M$.

If $A$ is the adjacency matrix of a strongly regular graph $G$, then Equation~(\ref{eqn:evals}) implies that $A^2$ is contained in the span of $\{I,A,J\}$. Since we also have $AJ = JA = kJ$, this further implies that any polynomial in $A$ is contained in this span. Letting $\bar{A} = J - I - A$ be the adjacency matrix of the complement of $G$, it is easy to see that this span is equal to the span of $\{I,A,\bar{A}\}$. Denoting by $E_\th$ and $E_\tau$ the projections onto the $\th$- and $\tau$-eigenspaces of $A$ respectively, we have that both of these projectors are contained in the span of $\{I,A,\bar{A}\}$. This means that $E_\th$ and $E_\tau$ have three distinct entries: those corresponding to vertices, edges, and non-edges of $G$.

The exact value of the entries of $E_\th$ and $E_\tau$ can be determined using a simple matrix identity. Specifically, 
an easy computation shows that the following holds for any real matrices $M$ and $N$ with the same dimensions:
\begin{equation}\label{eqn:identity}
\tr(M^T N) = \text{sum}(M \circ N).
\end{equation}
Note that we can drop the transpose when dealing with symmetric matrices, which will generally be the case for us.

We can now use Equation~(\ref{eqn:identity}) to compute the entries of $E_\tau$. For example, the entries of $E_\tau$ corresponding to edges are equal to
\[\frac{1}{nk}\text{sum}(A \circ E_\tau) = \frac{1}{nk}\tr(A E_\tau) = \frac{1}{nk}\tr(\tau E_\tau) = \frac{\tau m_\tau}{nk},\]
since the trace of a projector is equal to its rank. Similar computations for the other entries of $E_\tau$ reveal that
\[\left(E_\tau\right)_{uv} =
\begin{cases}
m_\tau/n & \text{if } u = v \\
\tau m_\tau/nk & \text{if } u \sim v \\
(-\tau - 1) m_\tau/n(n-k-1) & \text{if } u \not\simeq v
\end{cases}\]
One can also determine the entries of $E_\th$ in a similar manner, but we will not need this.

The proof of our main result makes use of of the projection $E_\tau$, but we will actually want to scale this matrix so that its diagonal entries are equal to one. Thus we define the \emph{cosine matrix} of a strongly regular graph $G$, denoted $E_G$, to be the matrix given as follows:
\[\left(E_G\right)_{uv} = \frac{n}{m_\tau} \left(E_\tau\right)_{uv} =
\begin{cases}
1 & \text{if } u = v \\
\tau/k & \text{if } u \sim v \\
(-\tau - 1)/(n-k-1) & \text{if } u \not\simeq v
\end{cases}\]

The key properties of $E_G$ that we will make use of are that it is positive semidefinite and that $(A - \tau I)E_G = 0$, both of which follow from the fact that it is a positive multiple of $E_\tau$.

Since the matrix $E_G$ is positive semidefinite with ones on the diagonal, it is the Gram matrix of some unit vectors that we can consider as being assigned to the vertices of the graph. The off diagonal entries of $E_G$ are then the cosines of the angles between these vectors, thus motivating the term ``cosine matrix". We refer to the values $\tau/k$ and $(-\tau - 1)/(n-k-1)$ as the \emph{adjacency and non-adjacency cosines} of a strongly regular graph, respectively. Note that for a primitive strongly regular graph $G$, its adjacency cosine is always contained in the interval $(-1,0)$, and its non-adjacency cosine is contained in the interval $(0,1)$. The latter follows from the fact, presented in the next section, that $n-k-1$ and $-\tau - 1$ are the largest and second largest eigenvalues of the complement of $G$ respectively.

Note that the parameters of a strongly regular graph determine its adjacency and non-adjacency cosines, but the converse is not true. Indeed, strongly regular graphs with parameter sets $(16, 10, 6, 6)$, $(26, 15, 8, 9)$, or $(36, 20, 10, 12)$ all have adjacency and non-adjacency cosines equal to $-1/5$ and $1/5$ respectively.

\subsection{Complements and some combinatorial properties}

It is easy to check that if $G$ is a strongly regular graph with parameters $(n,k,\l,\m)$, then the complement of $G$, denoted $\overline{G}$, is also a strongly regular graph with parameters $(n, \bk, \bl, \bm)$ where
\begin{align*}
\bk &= n - k - 1 \\
\bl &= n-2k-2+\m \\
\bm &= n - 2k + \l
\end{align*}
The eigenvalues of $\overline{G}$ are denoted by $\bk > \bth > \btau$. The latter two can be computed from the parameters of $\overline{G}$ using the identities in Section~\ref{subsec:algprops}, but it is easier to use the fact that the adjacency matrix of $\overline{G}$ is equal to $J - I - A$, where $A$ is the adjacency matrix of $G$. From this it follows that
\begin{align*}
\bth &= -\tau - 1 \\
\btau &= -\th - 1
\end{align*}


The last property of strongly regular graphs that we will need concerns the second neighborhoods of vertices. The \emph{second neighborhood} of a vertex $v$, denoted $N_2(v)$, is the set of vertices at distance exactly two from $v$. The following result is from~\cite{gardiner}, but we provide the proof for the reader's convenience:

\begin{lemma}\label{lem:N2connected}
Let $G$ be a primitive strongly regular graph. For any $v \in V(G)$, the subgraph of $G$ induced by $N_2(v)$ is connected.
\end{lemma}
\proof
Let $v \in V(G)$ and suppose that the subgraph induced by $N_2(v)$ is not connected. Let $C_1$ and $C_2$ be two distinct connected components of this induced subgraph, and let $u \in C_1$, $w \in C_2$. Since $G$ is strongly regular, both $u$ and $w$ share exactly $\mu$ common neighbors with $v$. But since $u$ and $w$ are also not adjacent, these are also the $\mu$ neighbors common to $u$ and $w$. It follows that all vertices of $N_2(v)$ share the same set $S$ of $\mu$ neighbors with $v$ and each other. Now suppose that $k > \mu$ and that $x \not\in S$ is a neighbor of $v$. Then $x$ shares the same number of neighbors with $v$ as any of the vertices in $S$, and thus is also adjacent to the same number of vertices of $N_2(v)$ as any vertex of $S$. But then we would have that $x$ is in $S$, a contradiction. Thus $k = \mu$ and $G$ is complete multipartite, a contradiction to the assumption that $G$ was primitive.\qeds

\section{Properties of Homomorphisms Between SRGs}\label{sec:main}

In this section we prove our main result that any homomorphism between strongly regular graphs with the same parameters is either an isomorphism or a coloring. However, we will first need to introduce the following construction:

\begin{definition}
Suppose that $M$ is a symmetric matrix with rows and columns indexed by some finite set $T$. For any set $S$ and function $\varphi: S \to T$, let $M^\varphi$ denote the matrix indexed by $S$ and defined entrywise as $\left(M^\varphi\right)_{uv} = M_{\varphi(u) \varphi(v)}$.
\end{definition}

It turns out that this construction preserves positive semidefiniteness:

\begin{lemma}\label{lem:psd}
Suppose $M$ is a positive semidefinite matrix indexed by some set $T$ and let $\varphi: S \to T$ for some set $S$. Then $M^\varphi$ is positive semidefinite.
\end{lemma}
\proof
Since $M$ is positive semidefinite, it is the Gram matrix of some multiset of vectors $\{p_w : w \in T\}$. In other words, $M_{ww'} = p_w^Tp_{w'}$. But then we have that $M^\varphi_{uv} = M_{\varphi(u)\varphi(v)} = p_{\varphi(u)}^Tp_{\varphi(v)}$. Thus $M^\varphi$ is the Gram matrix of the multiset of vectors $\{p_{\varphi(u)} : u \in S\}$, and is therefore positive semidefinite.\qeds

Using the above, we can prove the following which will be instrumental in proving our main result.

\begin{lemma}\label{lem:product}
Suppose $G$ and $H$ are strongly regular graphs with the same adjacency cosines. Let $A$ be the adjacency matrix of $G$ and $\tau$ its least eigenvalue. If $\varphi$ is a homomorphism from $G$ to $H$, then $(A - \tau I)E_H^\varphi = 0$.
\end{lemma}
\proof
First, recall that $(A- \tau I)E_G = 0$. Since $\varphi$ is a homomorphism and $G$ and $H$ have the same adjacency cosines, we have that $E_G$ and $E_H^\varphi$ agree on their diagonals and entries corresponding to the edges of $G$. Therefore,
\begin{align*}
\tr\left((A - \tau I)E_H^\varphi\right) &= \text{sum}\left((A - \tau I) \circ E_H^\varphi\right) \\
&= \text{sum}\left((A - \tau I) \circ E_G\right) \\
&= \tr\left((A - \tau I)E_G\right) = 0.
\end{align*}
Since both $A - \tau I$ and $E_H^\varphi$ are positive semidefinite (using Lemma~\ref{lem:psd} for the latter), the above implies that $(A - \tau I)E_H^\varphi = 0$.\qeds

Suppose that $G$ and $H$ are strongly regular graphs with equal adjacency cosines $\alpha$ and non-adjacency cosines $\beta$ and $\beta'$ respectively. If $\varphi$ is a homomorphism from $G$ to $H$, define the \emph{homomorphism matrix} of $\varphi$ to be $X := E_H^\varphi - E_G$. Then
\[X_{uv} = \begin{cases}
1 - \beta & \text{if } u \not\simeq  v \ \& \ \varphi(u) = \varphi(v) \\
\alpha - \beta & \text{if } u \not\simeq  v \ \& \ \varphi(u) \sim \varphi(v) \\
\beta' - \beta & \text{if } u \not\simeq  v \ \& \ \varphi(u) \not\simeq \varphi(v) \\
0 & \text{o.w.}
\end{cases}\]

Recall that $\alpha \in (-1,0)$ and $\beta, \beta' \in (0,1)$. Therefore we have that $1 - \beta > 0$ and $\alpha - \beta < 0$. The noteworthy property of the homomorphism matrix is that $(A - \tau I)X = 0$ where $A$ is the adjacency matrix of $G$ and $\tau$ its least eigenvalue. This follows immediately from the fact that $(A - \tau I)E_G = 0$ and $(A - \tau I)E_H^\varphi = 0$ by Lemma~\ref{lem:product}. The other important property of the homomorphism matrix is that it contains many zeros. This allows us to prove our main result:

\begin{theorem}\label{thm:main}
Let $G$ and $H$ be primitive strongly regular graphs with the same adjacency cosines equal to $\alpha$, and non-adjacency cosines equal to $\beta$ and $\beta'$ respectively. Suppose $\varphi$ is a homomorphism from $G$ to $H$. Then the following hold:
\begin{enumerate}
\item If $\beta > \beta'$, then $\varphi$ is a coloring.
\item If $\beta = \beta'$, then $\varphi$ is either a coloring or an isomorphism to an induced subgraph of $H$.
\end{enumerate}
\end{theorem}
\proof
Let $X$ be the homomorphism matrix of $\varphi$. Suppose that $\varphi$ is not a coloring. Then there exist vertices $u,v \in V(G)$ such that $\varphi(u) \not\simeq \varphi(v)$. Note that this implies that $u \not\simeq v$. For notational purposes, define the following sets:
\begin{align*}
C_1 &= \{w \in V(G): w \sim u, w \not\simeq v, \varphi(w) = \varphi(v)\} \\
C_2 &= \{w \in V(G): w \sim u, w \not\simeq v, \varphi(w) \sim \varphi(v)\} \\
C_3 &= \{w \in V(G): w \sim u, w \not\simeq v, \varphi(w) \not\simeq \varphi(v)\}
\end{align*}
Note that $C_1 \cup C_2 \cup C_3$ is the set of all neighbors of $u$ contained in $N_2(v)$.
 Since $\varphi$ is a homomorphism and we assumed that $\varphi(u) \not\simeq \varphi(v)$, we have that $C_1$ is empty. Now let $A$ be the adjacency matrix of $G$ and $\tau$ its least eigenvalue. Then $(A - \tau I)X = 0$ and therefore
\begin{align*}
0 = \left((A-\tau I)X\right)_{uv} &= \sum_{w \in V(G)} (A-\tau I)_{uw}X_{wv} \\
&= -\tau X_{uv} + \sum_{w \sim u} X_{wv} \\
&= -\tau (\beta' - \beta) + (1-\beta)|C_1| + (\alpha - \beta)|C_2| + (\beta' - \beta)|C_3| \\
&= -\tau (\beta' - \beta) + (\alpha - \beta)|C_2| + (\beta' - \beta)|C_3|.
\end{align*}
Now $\alpha - \beta < 0$, and $-\tau > 0$. Therefore, if $\beta > \beta'$ then every term above is non-positive, and the first term is strictly negative. This is a contradiction and so in this case no homomorphism that is not a coloring can exist. This proves the first claim.

If $\beta = \beta'$, then the above implies that $C_2$ is empty, and we already noted that $C_1$ is empty. Let us consider what this means. Since $C_1 \cup C_2 \cup C_3$ is the set of all neighbors of $u$ in $N_2(v)$, this implies that all such vertices $w$ satisfy $\varphi(w) \not\simeq \varphi(v)$. In other words, if $\varphi(u) \not\simeq \varphi(v)$, then $\varphi$ preserves non-adjacency between $v$ and every neighbor of $u$ in $N_2(v)$.

Now we can apply the above argument again, replacing $u$ with any neighbor of $u$ in $N_2(v)$. Since $N_2(v)$ is connected by Lemma~\ref{lem:N2connected}, iterating this argument implies that $\varphi$ preserves non-adjacency between $v$ and every vertex of $N_2(v)$. But now, for any $w \in N_2(v)$, we have that $\varphi(v) \not\simeq \varphi(w)$ and thus it must follow that $\varphi$ preserves non-adjacency between $w$ and every vertex of $N_2(w)$. Iterating again, and using the fact that $\overline{G}$ is connected, we see that $\varphi$ must preserve all non-adjacencies, i.e., it is an isomorphism to an induced subgraph of $H$.\qeds

As a corollary, we immediately obtain the following:

\begin{corollary}\label{cor:main}
If $G$ and $H$ are primitive strongly regular graphs with the same parameters, then any homomorphism from $G$ to $H$ is either a coloring or an isomorphism.
\end{corollary}
\proof
In this case we have that $\beta = \beta'$ in Theorem~\ref{thm:main}, and therefore any such homomorphism is a coloring or an isomorphism to an induced subgraph of $H$. However, since they have the same parameters, $G$ and $H$ have the same number of vertices. Therefore, any isomorphism to an induced subgraph of $H$ is simply an isomorphism to $H$.\qeds

Finally, we obtain a strengthening of the Cameron and Kazanidis conjecture:
\begin{corollary}
Every primitive strongly regular graph is a pseudocore.
\end{corollary}

\section{Cliques, Colorings, and the Homomorphism Order}\label{sec:homorder}

Since we now know that all homomorphisms between strongly regular graphs with the same parameters are either isomorphisms or colorings, it is worth considering the properties of the colorings. In order to distinguish them, we will refer to homomorphisms that are not also isomorphisms as \emph{proper homomorphisms}. We will see that, for a fixed parameter set, the proper homomorphisms between strongly regular graphs are not only required to be colorings, but colorings with a fixed number of colors.

To begin we will first need a well-known spectral bound on the size of a coclique (independent set) in a regular graph. This bound, known as the ``ratio bound" states that for a $k$-regular graph $G$ with $n$ vertices and least eigenvalue $\tau$, we have that
\[\alpha(G) \le \frac{n\tau}{\tau - k},\]
where $\alpha(G)$ denotes the maximum size of a coclique of $G$. Moreover, a coclique $S$ of $G$ meets this bound if and only if every vertex outside of $S$ has $-\tau$ neighbors in $S$. This bound was proven for strongly regular graphs by Delsarte~\cite{delsarte}, and extended to regular graphs by Hoffman~\cite{hoffman}. Cocliques meeting the bound are often referred to as \emph{Delsarte cocliques}, and cliques meeting the same bound for the complement are referred to as Delsarte cliques. We will present a very nice short proof of this bound due to Godsil (personal communication via the grapevine).

Let $A$ be the adjacency matrix of a $k$-regular, $n$-vertex graph $G$ with least eigenvalue $\tau$. Define the matrix
\[N = (A - \tau I) - \frac{k-\tau}{n}J.\]
By considering how it acts on an orthogonal basis of eigenvectors of $A$, it is not hard to show that the matrix $N$ is positive semidefinite. Therefore, $y^T N y \ge 0$ for any vector $y$. If $y$ is the characteristic vector of an independent set $S$ of $G$, then $y^T Ay = 0$. Therefore,
\[0 \le y^TNy = -\tau y^T y - \frac{k - \tau}{n}y^T J y = -\tau |S| - \frac{k-\tau}{n}|S|^2,\]
and the bound can be easily unraveled from here. Moreover, equality holds if and only if $y$ is a 0-eigenvector of $N$, from which the equality case condition can be deduced.

This bound on the independence number also provides a bound on the chromatic number. In particular, since $\chi(G) \ge n/\alpha(G)$ for any $n$-vertex graph $G$, if $G$ is $k$-regular with least eigenvalue $\tau$, then
\[\chi(G) \ge \frac{n}{\alpha(G)} \ge \frac{n}{n\tau/(\tau - k)} = 1 - \frac{k}{\tau}.\]
Colorings meeting this bound are referred to as Hoffman colorings, and such colorings (especially of strongly regular graphs) have received some attention in the literature~\cite{haemers96, haemers06}. We note here that the color classes in any Hoffman coloring must be Delsarte cocliques. This means that, in a Hoffman coloring, any vertex has $-\tau$ neighbors in each color class other than its own, in which it obviously has no neighbors. Therefore, the color classes of a Hoffman coloring form an \emph{equitable partition} of the graph. This is not directly relevant to what we will do here, but it is worth noting that Hoffman colorings appear to be quite special. Indeed, it is known that for a fixed $c \in \mathbb{N}$, only finitely many strongly regular graphs have Hoffman colorings with $c$ colors~\cite{haemersthesis}.

By taking complements, the ratio bound says that for a regular graph $G$, the maximum size of a clique in $G$ is at most $n\btau/(\btau - \bk)$, where $\bk$ and $\btau$ are the valency and least eigenvalue of $\overline{G}$. For strongly regular $G$, only minor arithmetical contortions are required to show that this is equal to $1 - k/\tau$, i.e.~the Hoffman bound on chromatic number. Therefore, for any strongly regular graph $G$, we have that
\begin{equation}\label{eqn:sandwich}
\omega(G) \le 1 - \frac{k}{\tau} \le \chi(G),
\end{equation}
where $\omega$ denotes the clique number. Importantly for us, this simultaneous bound on the clique and chromatic numbers of a strongly regular graph depends only on the parameters, not the specific graph. We are therefore able to prove the following:

\begin{lemma}\label{lem:properhoms}
Let $G$ and $H$ both be $SRG(n,k,\l,\m)$'s. There exists a proper homomorphism from $G$ to $H$ if and only if
\[\chi(G) = 1 - \frac{k}{\tau} = \omega(H),\]
i.e.~$G$ has a Hoffman coloring and $H$ contains a Delsarte clique.
\end{lemma}
\proof
Suppose there exists a proper homomorphism from $G$ to $H$. By Corollary~\ref{cor:main}, this homomorphism must be a coloring. Therefore, using Equation~(\ref{eqn:sandwich}), we have that
\[1 - \frac{k}{\tau} \le \chi(G) \le \omega(H) \le 1 - \frac{k}{\tau}.\]
The converse is trivial.\qeds

Note that the above lemma implies that if $G$ and $H$ are non-isomorphic $SRG(n,k,\l,\m)$'s, then $G \to H$ if and only if $\chi(G) = 1 - k/\tau = \omega(H)$. We also obtain the following corollary giving an if and only if condition for when a strongly regular graph is a core:

\begin{corollary}
If $G$ is a strongly regular graph, then $G$ is NOT a core if and only if
\[\omega(G) = 1 - \frac{k}{\tau} = \chi(G).\]
In this case the core of $G$ is a complete graph of size $1 - \frac{k}{\tau}$.
\end{corollary}

\subsection{Types and the homomorphism order}

The result of Lemma~\ref{lem:properhoms} suggests a useful partition of strongly regular graphs of a fixed parameter set. Namely, to classify them according to which subset of $\{\omega(G), \chi(G)\}$ meet the Hoffman bound. We therefore propose the following four ``types" of strongly regular graphs:
\begin{itemize}
\item Type A: $\omega(G) < 1 - \frac{k}{\tau} = \chi(G)$;
\item Type B: $\omega(G) = 1 - \frac{k}{\tau} = \chi(G)$;
\item Type C: $\omega(G) = 1 - \frac{k}{\tau} < \chi(G)$;
\item Type X: $\omega(G) < 1 - \frac{k}{\tau} < \chi(G)$.
\end{itemize}
The existence of a homomorphism between any two non-isomorphic $SRG(n,k,\l,\m)$'s is determined by their types: Any graph of type A or B has homomorphisms to any graph of type B or C. There are no other homomorphisms between non-isomorphic $SRG(n,k,\l,\m)$'s. Furthermore, all graphs of type A, C, or X are cores, and all graphs of type B have complete graphs of size $1 - k/\tau$ as their cores. Summarizing these observations, we have the following Hasse diagram of the homomorphism order of $SRG(n,k,\l,\m)$'s:

\begin{center}
\begin{figure}[h]
\includegraphics{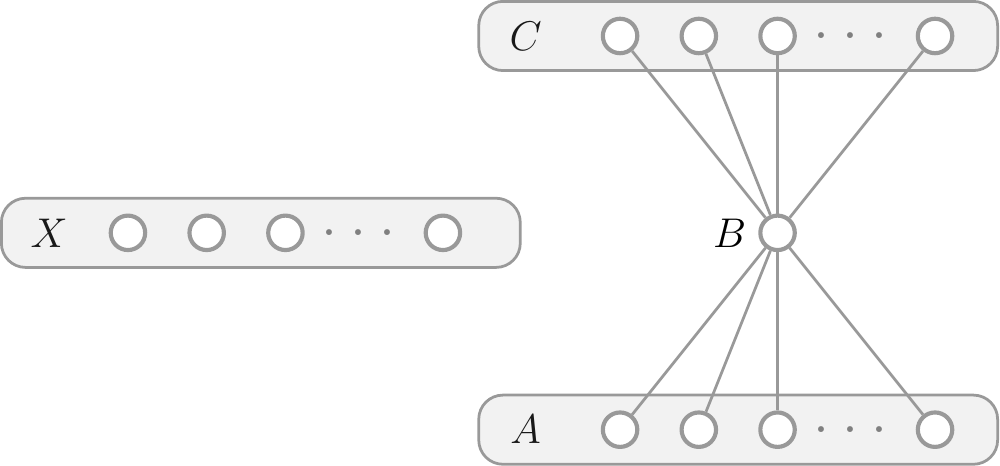}
\caption{Homomorphism order of $SRG(n,k,\l,\m)$'s.}
\end{figure}
\end{center}
Note that the type B graphs are represented by a single node in the above diagram since they are all homomorphically equivalent, whereas graphs of any other fixed type are incomparable (have no homomorphisms in either direction between them).

The four types defined above can also be defined purely in terms of the homomorphisms among $SRG(n,k,\l,\m)$'s, without explicitly referring to the Hoffman bound or clique or chromatic number. The type B $SRG(n,k,\l,\m)$'s are those which are homomorphically equivalent to at least one other $SRG(n,k,\l,\m)$. Of the remaining graphs, those of type A are the ones with homomorphisms to the type B graphs, those of type C admit homomorphisms from the type B graphs, and the type X graphs are incomparable to all other $SRG(n,k,\l,\m)$'s. The problem with this definition is that it assumes that graphs of type B exist for any given parameter set. For instance, if all of the SRG's for a parameter set are incomparable, then there may exist only type A graphs for this parameter set, or both type A and type X graphs (or some other combination of types).

If the Hoffman bound is not an integer, then neither the clique nor chromatic number can meet this bound with equality, and therefore only graphs of type X can occur. This happens for conference graphs of non-square order, since these have $\tau$ equal to an irrational number. However, this can also occur for other parameter sets. Some examples include $(10,3,0,1)$, $(16,5,0,2)$, $(21,10,3,6)$, $(26,10,3,4)$, $(36,14,4,6)$, and $(36,21,10,15)$, for all of which there do exist strongly regular graphs. Also note that if the Hoffman bound of the complementary parameter set is not an integer, then there can be no Delsarte cocliques, and therefore no Hoffman colorings. Therefore, for such parameter sets, there will only be type C and/or X graphs.

Computations reveal that there are parameter sets which contain only graphs of a single type. Examples of this for each type, including an example for type X where the Hoffman bound is an integer, are given below:
\begin{itemize}
\item Type A - $(27,16,10,8)$;
\item Type B - $(49,12,5,2)$;
\item Type C - $(45,32,22,24)$;
\item Type X - $(16,10,6,6)$.
\end{itemize}

On the other hand, there are also parameter sets having all four types. Some examples include $(36,20,10,12), (45,12,3,3)$, and $(64,18,2,6)$. In general, for the strongly regular graphs we performed computations on, which were obtained from Ted Spence's webpage~\cite{spence}, almost all of them were either type C or X. This seems to indicate that having a Hoffman coloring is a rare property for a strongly regular graph, but having a Delsarte clique is not. The latter observation is perhaps not so surprising since it is known that all strongly regular graphs arising as point graphs of partial geometries have Delsarte cliques.

The computations for the above were done in Sage~\cite{sage}. One only needs to determine if the given strongly regular graph has a clique of a certain size and/or coloring with certain number of colors. For the former, the built in clique number routine is very fast, and so there is no problem finding the clique number of all the strongly regular graphs from Ted's webpage. This is not the case for chromatic number. Sage's built in coloring routines seem to be far too slow to be of any use for this endeavor. However, there is a GAP package called Digraphs~\cite{digraphs} developed by researchers at The University of St Andrews, and the coloring routine in this package works very quickly in comparison. In fact, it is hard to overstate how much faster it seems to be.

\section{Vector Colorings and the Lov\'{a}sz $\vartheta$ Function}\label{sec:vectcols}

In this section we will see that some of the results of Section~\ref{sec:main} are part of a more general theory involving semidefinite programs and the Lov\'{a}sz theta number of a graph. 

For a graph $G$ and real number $t \ge 2$, a \emph{strict vector $t$-coloring} of $G$ is an assignment, $u \mapsto p_u$, of unit vectors to the vertices of $G$ such that
\[p_u^T p_v = \frac{-1}{t-1} \text{ for all } u \sim v.\]
If we drop the ``strict", then we only require that the inner product above is upper bounded by the righthand side. We note however that for strongly regular graphs, every optimal vector coloring is also a strict vector coloring~\cite{UVC1}. For a non-empty graph $G$, its \emph{strict vector chromatic number} is the minimum $t \ge 2$ such that $G$ admits a strict vector $t$-coloring. For empty graphs, this parameter is defined to be equal to 1. The strict vector chromatic number was defined by Karger, Motwani, and Sudan~\cite{KMS}, and they showed that it is equal to the Lov\'{a}sz theta number of the complement graph. The Lov\`{a}sz theta number is typically denoted by $\vartheta$, and so we will use $\bartheta(G) := \vartheta(\overline{G})$ to denote the strict vector chromatic number of $G$. We will give two of the more well known formulations of the Lov\'{a}sz theta number in Section~\ref{subsec:sdp}.

By considering the Gram matrix of vectors in a strict vector coloring, it is easy to see that $G$ has a strict vector $t$-coloring if and only if there exists a positive semidefinite matrix $M$ indexed by the vertices of $G$ such that
\[M_{uv} =
\begin{cases}
1 & \text{if } u = v \\
\frac{-1}{t-1} & \text{if } u \sim v
\end{cases}\]
Using this interpretation, it is not difficult to see that a complete graph on $n$ vertices has strict vector chromatic number equal to $n$. It is also now apparent that the matrices $E_G$ and $E_H^\varphi$ from Section~\ref{sec:main} were Gram matrices of strict vector colorings.

Suppose that $G$ and $H$ are graphs and that $w \mapsto p_w$ for $w \in V(H)$ is a strict vector $t$-coloring of $H$. If $\varphi$ is a homomorphism from $G$ to $H$, then it is easy to see that $u \mapsto p_{\varphi(u)}$ for $u \in V(G)$ is a strict vector $t$-coloring of $G$ (note that this is the exact construction used in the proof of Lemma~\ref{lem:psd} to show that $M^\varphi$ is positive semidefinite). It follows that if $G \to H$, then $\bartheta(G) \le \bartheta(H)$, i.e.~the strict vector chromatic number is \emph{homomorphism monotone}. In particular, using the fact that $\bartheta(K_n) = n$, this implies the well known ``sandwich theorem'':
\[\omega(G) \le \bartheta(G) \le \chi(G).\]

\subsection{Semidefinite programming}\label{subsec:sdp}

One of the many useful properties of the Lov\'{a}sz theta number is that it can be written as a semidefinite program that satisfies strong duality. This provides us with both a minimization and maximization program for this parameter:

\[
\begin{array}[t]{lcc}
 & \text{PRIMAL} & \text{DUAL} \vspace{.05in} \\
\bartheta(G) \ = & \begin{array}[t]{ll}
\min & t \\
\text{s.t.} & M_{uu} = t-1 \text{ for } u \in V(G) \\
 & M_{uv} = -1 \text{ for } u \sim v \\
 & M \succeq 0
\end{array} = & \begin{array}[t]{ll}
\max & \text{sum}(B) \\
\text{s.t.} & B_{uv} = 0 \text{ for } u \not\simeq v \\
 & \tr(B) = 1 \\
 & B \succeq 0
\end{array}
\end{array}
\]

Note that a feasible solution of value $t$ for the primal program above is exactly $(t-1)$ times the Gram matrix of a strict vector $t$-coloring of $G$, and so we see that these are equivalent definitions of $\bartheta$.

Suppose that $M$ and $B$ are feasible solutions to the above primal and dual formulations of $\bartheta$ with objective values $P$ and $D$ respectively. Then,
\[\tr(MB) = \text{sum}(M \circ B) = (P-1)\tr(B) - \left[\text{sum}(B) - \tr(B)\right] = P-D.\]
It thus follows that if $M$ and $B$ are feasible solutions for the primal and dual programs respectively, then they are both optimal if and only if $\tr(MB) = 0$ if and only if $MB = 0$. This is in fact just the complementary slackness condition for these semidefinite programs.

For any graph $G$ with adjacency matrix $A$ and least eigenvalue $\tau$, the matrix $A - \tau I$ meets the first and third conditions for the dual program above. If we let $B$ be the positive scaling of $A - \tau I$ that has trace one, then $B$ is a feasible solution to the dual. If $G$ is strongly regular, then we have seen in Section~\ref{subsec:projections} that the cosine matrix of $G$, $E_G$, is constant on the diagonal, and is a negative constant on entries corresponding to edges of $G$. Therefore, up to a scalar multiple, this is a feasible solution to the primal program for $\bartheta(G)$. If we let $M$ denote this scalar multiple of $E_G$, then it is obvious that $MB = 0$. Therefore these are both optimal solutions to their respective programs. It is then only a matter of arithmetic to show that $\bartheta(G)$ is equal to our old friend the Hoffman bound for any strongly regular graph $G$.

We can now see Lemma~\ref{lem:product} for what it is:\footnote{All instances of the phrase ``up to a scalar'' have been removed from the following so that the printers do not run out of ink.} The strongly regular graph $G$ has feasible solutions $E_G$ and $A - \tau I$ to the primal and dual respectively, and these must be optimal since they multiply to 0. Similarly, the cosine matrix $E_H$ is an optimal primal solution for $H$, and $E_H^\varphi$ is the Gram matrix of the strict vector coloring of $G$ obtain by composing $\varphi$ with the strict vector coloring of $H$ whose Gram matrix is $E_H$. Since both graphs are strongly regular with the same parameters, they have the same strict vector chromatic number and therefore $E_H^\varphi$ is an optimal primal solution for $G$. Finally, since $A - \tau I$ was already shown to be an optimal dual solution for $G$, we have that $(A - \tau I)E_H^\varphi = 0$.

Of course, a similar technique can be applied to any homomorphism between graphs with the same strict vector chromatic number. But the primal and dual solutions for the two graphs will likely not be as nice as in the strongly regular case. The key feature of the primal solutions we used is that their entries depend only on whether the corresponding vertices are equal, adjacent, or non-adjacent. Most graphs will not have an optimal primal solution of this form.

On the other hand, distance regular graphs also have $E_\tau$ and $A - \tau I$ as optimal primal and dual solutions, and the $uv$-entry of the matrix $E_\tau$ only depends on the distance between vertices $u$ and $v$. Thus, distance regular graphs are a natural choice for attempting to generalize our main theorem. Indeed, strongly regular graphs are exactly distance regular graphs of diameter two. However, the analysis seems more difficult in this case, since the matrix $E_H^\varphi - E_G$ will potentially have a different nonzero entry for every way in which the homomorphism $\varphi$ can change the distance between two vertices. This is actually the same for our case, but for us there were only two such possibilities.

Another possible route for generalization would be to consider \emph{directed} strongly regular graphs. These were introduced in~\cite{duval} and have been given a fair amount of attention in the literature. Since homomorphisms extend naturally to directed graphs, and many of the algebraic properties of strongly regular graphs have analogs in the directed case~\cite{GHM}, it seems plausible that our main result could be generalized to this larger class of graphs.


\section{A Generalization}\label{sec:generalization}

We did not make extensive use of the fact that $H$ was a strongly regular graph in the proof of our main result, nor the lemmas leading up to it. If we let $G$ be an $SRG(n,k,\l,\m)$, then the only thing we required of $H$ in our arguments is that the matrix $I + \alpha A_H + \beta' \overline{A}_H$, where $\alpha$ is the adjacency cosine of $G$ and $\beta'$ is at most the non-adjacency cosine of $G$, is positive semidefinite. The proof of the main result now proceeds exactly as before.

The assumption that $I + \alpha A_H + \beta' \overline{A}_H$ is positive semidefinite implies that $H$ admits a strict vector coloring of value $1 - 1/\alpha = \bartheta(G)$. Since we also assumed that $G  \to H$, this must be an optimal strict vector coloring of $H$. This inspires the following definition. For real numbers $\alpha$ and $\beta$, we say that $H$ is an $(\alpha,\beta)$-graph if $I + \alpha A_H + \beta \overline{A}_H$ is the Gram matrix of an optimal strict vector coloring of $H$. Note that this implies that $\alpha \in [-1,0)$. We can now succinctly state the above discussed generalization of our main result: 

\begin{theorem}\label{thm:generalization}
Suppose $G$ is an strongly regular graph with adjacency and non-adjacency cosines $\alpha$ and $\beta$ respectively, and that $H$ is an $(\alpha,\beta')$-graph. Let $\varphi$ be a homomorphism from $G$ to $H$. Then the following hold:
\begin{enumerate}
\item If $\beta > \beta'$, then $\varphi$ is a coloring.
\item If $\beta = \beta'$, then $\varphi$ is either a coloring or an isomorphism to an induced subgraph of $H$.
\end{enumerate}
\end{theorem}

Note that in the case of a coloring, the image of $\varphi$ must be a maximum clique of $H$ of size $1 - \frac{k}{\tau}$. In either case, the image of $\varphi$ must have strict vector chromatic number equal to that of both $G$ and $H$, namely $1 - \frac{k}{\tau}$.


\section{Discussion}

The main purpose of this work was to prove the conjecture of Cameron \& Kazanidis. However, our results have several other implications and raise certain questions. We will discuss some of these here.

Since all but finitely many strongly regular graphs with fixed least eigenvalue are the point graphs are partial geometries, these geometric graphs warrant some consideration with respect to our results. We mentioned previously that geometric graphs always have Delsarte cliques. This is because the Hoffman bound for these graphs is equal to the size of a line in the underlying partial geometry, and thus the points on a line induce a Delsarte clique, though there may be others. It follows from this that all geometric graphs are of types B or C. Therefore, a geometric graph is type B if and only if it has a Hoffman coloring, and otherwise is type C. Recall that every color class in a Hoffman coloring is a Delsarte coclique. For geometric graphs, it is known that a Delsarte coclique corresponds to a set of points in the underlying partial geometry that meets every line exactly once, and vice versa. Such an object is called an \emph{ovoid}. Therefore, a Hoffman coloring of a geometric graph is a partition of its partial geometry into ovoids. A partition into ovoids is, for obvious reasons\footnote{To geometers, presumably.}, called a \emph{fan}. So we see that the point graph of a partial geometry is type B if and only if the geometry has a fan, and otherwise the graph is type C.

In light of the generalization of our main result presented in Section~\ref{sec:generalization}, it is interesting to ask what graphs are $(\alpha, \beta)$-graphs for which real numbers $\alpha$ and $\beta$. We are presently preparing a paper addressing this question, but we will discuss some basic points here. First, we have seen that strongly regular graphs are $(\alpha, \beta)$-graphs for $\alpha = \tau/k$ and $\beta = \bth/\bk$. As we mentioned in Section~\ref{subsec:projections} it is possible for different parameter sets to result in the same values of both $\tau/k$ and $\bth/\bk$. This brings us to an interesting question: for fixed $\alpha$ and $\beta$, are there an infinite number of $(\alpha, \beta)$-graphs? If we restrict to strongly regular graphs, it turns out the answer is no. This is because, as we show in our upcoming paper, the second largest eigenvalue of a regular $(\alpha,\beta)$-graph is determined by $\alpha$ and $\beta$. Thus the least eigenvalue of its complement is determined. So for fixed $\alpha$ and $\beta$, the least eigenvalue of the complement of a strongly regular $(\alpha,\beta)$-graph is fixed, and thus Neumaier's result can be applied. One can then simply check the infinite families to see that these do not provide infinitely many $(\alpha, \beta)$-graphs.

In the positive direction, any graph which is transitive on its non-edges is an $(\alpha, \beta)$-graph for some values of $\alpha$ and $\beta$. This is because the Gram matrix of any optimal strict vector coloring of a non-edge-transitive graph can be ``smoothed out" on the non-edges by taking a uniform convex combination of the Gram matrix conjugated by permutation matrices representing automorphisms of the graph. This provides a large class of $(\alpha, \beta)$-graphs that includes many graphs which are not strongly regular.


The fact that every strongly regular graph is a pseudocore has implications in the study of synchronizing groups. A permutation group $\Gamma$ acting on a set $S$ \emph{synchronizes} a function $f$ from $S$ to itself if the monoid generated by $\Gamma$ and $f$ contains a transformation whose image is a single element of $S$. The group $\Gamma$ is said to be \emph{synchronizing} if it synchronizes every function that is not a permutation. This definition is motivated by concerns in the theory of finite automata, in particular the \v{C}ern\'{y} conjecture. In~\cite{synch}, Cameron et. al.~define \emph{almost synchronizing} permutation groups as those which synchronize all functions which are non-uniform, i.e.~whose preimages are not all the same size. They note that the automorphism group of any vertex transitive pseudocore is almost synchronizing whenever it is primitive. Therefore, our main result shows that the automorphism group of any vertex transitive strongly regular graph is almost synchronizing whenever it is primitive. In particular, they note\footnote{Cameron et. al.~received a preprint of this manuscript before it became publicly available.} that this implies any primitive group with permutation rank 3 is almost synchronizing.

In~\cite{symcores}, the \emph{hull} of a graph was introduced by Cameron \& Kazanidis in order to prove that rank 3 graphs are core-complete. The hull of a graph $G$ has the same vertex set as $G$, and two vertices are adjacent in the hull if there does not exist an endomorphism which identifies these vertices. In particular, this means that every edge of $G$ is an edge of its hull. Cameron \& Kazanidis proved several results about the hull of a graph, showing that it is in some sense a dual notion to that of the core. It therefore may be natural to ask whether the hull of a strongly regular graph is always either the graph itself or a complete graph. This turns out to not be the case, and in fact we have found through direct computations that there are strongly regular graphs whose hulls are not even regular. We will not present a specific case, but we note that there are examples among the 23 type B $SRG(45,12,3,3)$'s.

\paragraph{Acknowledgements:} I would like to thank Chris Godsil, Brendan Rooney, Robert \v{S}\'{a}mal, and Antonis Varvitsiotis for all I learned through our work on vector colorings. In particular, I thank Robert for first showing me a paper about unique vector colorings which initiated this collaboration. Brendan was the unfortunate soul that first checked the proof of the main result for correctness, and I thank him for that. Chris, Krystal Guo, and Gordon Royle also read early versions of this work and I thank them for their helpful input. I would also like to thank Laura Man\v{c}inska for encouraging me to keep going after I found a mistake in my original proof.




\bibliographystyle{plainurl}

\bibliography{Homomorphisms_of_SRGs.bbl}

\end{document}